\documentclass[12pt,a4paper]{amsart}
\usepackage[abbrev]{amsrefs}
\usepackage{amssymb,amsmath}

\theoremstyle{plain}
  \newtheorem{thm}{Theorem}[section]
  \newtheorem{mthm}[thm]{Main Theorem}
  \newtheorem{lem}[thm]{Lemma}
  
  \newtheorem{cor}[thm]{Corollary}
  \newtheorem{prop}[thm]{Proposition}

\theoremstyle{definition}
  \newtheorem{defn}[thm]{Definition}
  
  \newtheorem{ex}[thm]{Example}

  \newtheorem{ques}[thm]{Question}

\theoremstyle{remark}
  \newtheorem{rem}[thm]{Remark}

\numberwithin{equation}{section}

\DeclareMathOperator{\diam}{diam}

\DeclareMathOperator{\dGH}{{\it d}_{{\rm GH}}}
\DeclareMathOperator{\dH}{{\it d}_{{\rm H}}}

\DeclareMathOperator{\CAT}{CAT}
\DeclareMathOperator{\rad}{rad}

\newcommand{\RP}{\R P}

\newcommand{\cG}{\mathcal{G}}
\newcommand{\cH}{\mathcal{H}}

\newcommand{\cM}{\mathcal{M}}

\newcommand{\cP}{\mathcal{P}}

\newcommand{\field}[1]{\mathbb{#1}}

\newcommand{\R}{\field{R}}

\begin{document}

\title[Compactification of the GH space]
{A natural compactification\\
of the Gromov-Hausdorff space}

\thanks{This work was supported by JSPS KAKENHI Grant Number 19K03459 and 19J10866.}

\author{Hiroki Nakajima}
\address{Institute for Excellence in Higher Education, Tohoku University,
Sendai 980-8576, Japan}
\email{hiroki.nakajima.a1@tohoku.ac.jp}

\author{Takashi Shioya}
\address{Mathematical Institute, Tohoku University, Sendai 980-8578, JAPAN}
\email{shioya@math.tohoku.ac.jp}

\begin{abstract}
In this paper, we introduce a pseudometric
on the family of isometry classes of (extended) metric spaces.
Using it, we obtain a natural compactification
of the Gromov-Hausdorff space, which is compatible with ultralimit.
\end{abstract}

\date{\today}



\maketitle

\section{Introduction}

The Gromov-Hausdorff (GH for short) distance is one of the most
natural concepts of distance between two compact metric spaces.
However, it is not effective for noncompact metric spaces,
and there are many interesting examples of sequence of metric spaces
that are not (even pointed) GH convergent.
As an effective generalization of GH limit,
we have the concept of ultralimit.
However, ultralimit does not induce a topology on
the space of metric spaces,
which is a disadvantage of ultralimit.
In this paper, we introduce a pseudometric
on the family of isometry classes of extended metric spaces.
Using it, we obtain a natural compactification
of the GH space, where the GH space is the the family of
isometry classes of compact metric spaces
equipped with the GH metric.
Since our compactification has a metric,
it is second-countable.

Let us mention our results more precisely.
A metric space whose metric takes values in
$[\,0,+\infty\,]$ is called an \emph{extended metric space}.
The reason why we consider extended metric spaces
is that the limit of metric spaces is an extended metric space
in general.

\begin{defn}[Relation $\precsim$]
For two extended metric spaces $X$ and $Y$,
we say $X \precsim Y$ if
for any finite sequence of points
$x_1,\dots,x_N \in X$,
there exists a sequence $y_{i,k} \in Y$, $i=1,2,\dots,N$, $k=1,2,\dots$,
such that
\[
d_X(x_i,x_j) \le \liminf_{k\to\infty} d_Y(y_{i,k},y_{j,k}).
\]
for any $i,j = 1,\dots,N$.
We say that $X$ is \emph{equivalent} to $Y$ and write $X \sim Y$
if $X \precsim Y$ and $Y \precsim X$ both hold.
\end{defn}

For example, if $X$ is a subspace of $Y$, then $X \precsim Y$.

Let $\cM$ denote the family of isometry classes
of extended metric spaces.
The relation $\precsim$ is a preorder relation on $\cM$,
and $\sim$ an equivalence relation on $\cM$.
In the case where $X$ and $Y$ are compact metric spaces,
we have $X \sim Y$ if and only if $X$ and $Y$ are isometric
to each other.
Any extended metric space is equivalent to its completion
(also to any ultrapower).

In this paper, we construct a pseudometric, denoted by $\rho$,
in Section \ref{sec:rho}.
Since $X \sim Y$ $\Leftrightarrow$ $\rho(X,Y) = 0$
(Corollary \ref{cor:rho-sim}),
the $\rho$ induces a metric on the quotient space
$\cM/\!\sim$.
One of our main theorems is stated as follows.

\begin{mthm} \label{mthm:main}
\begin{enumerate}
\item $(\cM/\!\sim,\rho)$ is a compact metric space.
\item The natural projection from the GH space to $\cM/\!\sim$
is a $3$-Lipschitz topological embedding map.
\item The projection in {\rm(2)} has dense image in
$(\cM/\!\sim,\rho)$.
\end{enumerate}
\begin{enumerate}
\item[(4)] For any $\rho$-convergent sequence
of extended metric spaces,
any ultralimit of it is equivalent to the $\rho$-limit.
\end{enumerate}
\end{mthm}

(1)---(3) of the theorem imply that
$(\cM/\!\sim,\rho)$ is a compactification of the GH space.

Different from ultralimit, it is sometimes easy to handle
the convergence of concrete sequences of metric spaces
and their limits in $\cM/\!\sim$ with respect to $\rho$.
In Section \ref{sec:examples}, we give several examples
of convergent sequences, such as, sequences of manifolds with
dimensions tending to infinity.

Note that $\cM$ is not a set,
however we have a representation of $\cM/\!\sim$ defining
it as a set (see Remark \ref{rem:cM} for the detail).

(2) says that
$\rho(X,Y) \le 3\dGH(X,Y)$ for any compact metric spaces
$X$ and $Y$.
We remark that the converse inequality like
$\dGH(X,Y) \le C\,\rho(X,Y)$
cannot be expected, where $C$ is a constant.
In fact, this kind of inequality does not hold for any compactification
of the GH space (see Remark \ref{rem:dGH-rho-2}).
By the same reason, if we extend the GH metric to
a (pseudo)metric on the family of metric spaces, then
it is impossible to obtain a compactification.

We also remark that our metric $\rho$ does not fit
to study (non-extended) metric spaces with infinite diameter,
because such metric spaces are all equivalent to each other
(see Proposition \ref{prop:S-infty-infty}).
To overcome this problem, we introduce a pointed version
of $\rho$.  We define a relation $\sim$ between
pointed metric spaces in the same manner
(see Definition \ref{defn:pointed-sim}), and then consider
the following.

\begin{defn}[Strongly equivalence]
For two pointed metric space $(X,x_0)$ and $(Y,y_0)$,
we say that $(X,x_0)$ is \emph{strongly equivalent} to $(Y,y_0)$
and write $(X,x_0) \simeq (Y,y_0)$
if $(B(x_0,r),x_0) \sim (B(y_0,r),x_0)$ for any $r > 0$,
where $B(x_0,r)$ is the closed metric ball centered at $x_0$
and of radius $r$.
\end{defn}

Denote by $\cG_0$ the family of isometry classes of
pointed geodesic metric spaces,
and by $\cP\cG_0$ the family of isometry classes of
pointed proper geodesic metric spaces
equipped with the pointed GH topology.
In Section \ref{sec:sigma}, we introduce a pseudometric,
denoted by $\rho_0$, on $\cG_0$
which satisfies that $\rho_0((X,x_0),(Y,y_0)) = 0$ $\Leftrightarrow$
$(X,x_0) \simeq (Y,y_0)$ (Theorem \ref{thm:sigma}),
and then prove the following.

\begin{mthm} \label{mthm:main-pointed}
\begin{enumerate}
\item $(\cG_0/\!\sim,\rho_0)$ is a compact metric space.
\item A natural projection from $\cP\cG_0$ to $\cG_0/\!\simeq$
is topological embedding map.
\item[(3)] For any $\rho_0$-convergent sequence
of pointed geodesic metric spaces,
any pointed ultralimit is strongly equivalent to the $\rho_0$-limit.
\end{enumerate}
\end{mthm}

By (1) and (2) of Main Theorem \ref{mthm:main-pointed},
the closure of the image of the embedding map in (2)
is a compactification of $\cP\cG_0$.

Different from $\rho$, we see that
$\rho_0((X,x_0),(Y,y_0)) > 0$ in general
even if $X$ and $Y$ have infinite diameter.

The construction of $\rho$ is as follows.
Consider the set, say $\Pi$, of subsets of the GH space satisfying
suitable conditions, and obtain a bijection between $\Pi$ and
$\cM/\!\sim$ in a natural way.  We then define $\rho$ by using $\Pi$.
Gromov \cite{Gmv:green} (see also \cites{Sy:mmg, Sy:mmlim})
introduced a compactification of the space of
metric measure spaces with the concentration topology.
The idea of $\Pi$ came from Gromov's compactification,
but the detail of our formulation is much different from Gromov's.
Our essential idea of the definition of $\rho(X,Y)$ is to compare
all finite subsets of $X$ and those of $Y$,
which does not appear in Gromov's theory.

We extend $\rho$ to a pointed version
to define $\rho_0$ as the integral of $\rho$-distance
between pointed metric balls with a certain weight function.
Applying a pointed version of Main Theorem \ref{mthm:main},
we prove Main Theorem \ref{mthm:main-pointed}.

\section{Preliminaries}

We refer to \cite{BBI} for the basics of metric geometry.
In this paper, $\dGH$ denotes the GH metric on the GH space.
For a metric space $X$, we denote its metric by $d_X$.
For extended real numbers $a$ and $b$,
we denote by $a \wedge b$ (resp.~$a \vee b$)
the minimum (resp.~maximum) of $a$ and $b$.

\begin{defn}[Weak Hausdorff convergence, \cite{Gmv:green}]
Let $X$ be a metric space and
$S,S_n \subset X$, $n=1,2,\dots$, closed subsets.
We say that $S_n$ converges to $S$ in the \emph{weak Hausdorff}
sense if the following (i) and (ii) are satisfied.
\begin{enumerate}
\item[(i)] $\lim_{n\to\infty} d_X(x,S_n) = 0$ for any $x \in S$.
\item[(ii)] $\liminf_{n\to\infty} d_X(x,S_n) > 0$
for any $x \in X \setminus S$.
\end{enumerate}
\end{defn}

We here set $d_X(x,S_n) := \inf_{x' \in S_n} d_X(x,x')$.
Note that $S_n$ and $S$ may be empty and
we set $d_X(x,\emptyset) = \infty$.

\begin{lem}[\cite{Sy:mmg}*{Lemma 6.6}]\label{lem:weakH-cpt}
The set of closed subsets of a metric space
is sequentially compact with respect to
the weak Hausdorff convergence.
\end{lem}

A standard discussion proves the following lemma.
We omit the proof.

\begin{lem} \label{lem:weak-H}
For given compact subsets $S_n, S \subset X$, $n=1,2,\dots$,
the following {\rm(1)} and {\rm(2)} are equivalent to each other.
\begin{enumerate}
\item $S_n$ converges to $S$ in the Hausdorff sense as $n\to\infty$.
\item $S_n$ converges to $S$ in the weak Hausdorff sense as $n\to\infty$.
\end{enumerate}
\end{lem}

\begin{defn}[$1$-Lipschitz up to $\varepsilon$]
Let $\varepsilon \ge 0$ be a real number.
A map $f : X \to Y$ between two metric spaces $X$ and $Y$
is said to be \emph{$1$-Lipschitz up to $\varepsilon$}
if
\[
d_Y(f(x),f(x')) \le d_X(x,x') + \varepsilon
\]
for any $x,x' \in X$.
\end{defn}

Denote by $\|\cdot\|_\infty$ the $l^\infty$-norm on $\R^N$.

\begin{lem}[see the proof of \cite{Sy:mmg}*{Lemma 5.4}]
\label{lem:f-tilde}
If a map $f : X \to (\R^N,\|\cdot\|_\infty)$ is
$1$-Lipschitz up to a real number $\varepsilon \ge 0$, then
there exists a $1$-Lipschitz map
$\tilde{f} : X \to (\R^N,\|\cdot\|_\infty)$ such that
$\|\tilde{f}(x)-f(x)\|_\infty \le \varepsilon$
for any $x,x' \in X$.
\end{lem}

\section{Construction of metric $\rho$} \label{sec:rho}

In this section, we define the metric $\rho$
and prove Main Theorem \ref{mthm:main}.

\begin{lem} \label{lem:precsim-cpt}
Let $X$ and $Y$ be two compact metric spaces
with $X \precsim Y$.
Then, there exists a $1$-Lipschitz surjective map
from a closed subset of $Y$ to $X$.
\end{lem}

\begin{proof}
We find a dense countable subset $\{x_i\}_{i=1}^\infty \subset X$.
It follows from $X \precsim Y$ and the compactness of $Y$
that for any number $N$ there is a sequence $\{y_i^N\}_{i=1}^N \subset Y$
such that $d_X(x_i,x_j) \le d_Y(y_i^N,y_j^N)$ for any $i,j = 1,2,\dots,N$.
Replacing $\{N\}$ with a subsequence, we assume that
for each $i$, the point $y_i^N$ converges to a point, say $y_i \in Y$,
as $N\to\infty$.  We have $d_X(x_i,x_j) \le d_Y(y_i,y_j)$ for any $i,j$.
Setting $f(y_i) := x_i$, we see the $1$-Lipschitz continuity of $f$.
Extend $f$ to a $1$-Lipschitz map from the closure of $\{y_i\}$.
Then, $f$ is surjective.
This completes the proof.
\end{proof}

\begin{lem} \label{lem:prec-isom}
Let $X$ and $Y$ be two compact metric spaces.
If $X \sim Y$, then $X$ and $Y$ are isometric to each other.
\end{lem}

\begin{proof}
Assume $X \sim Y$ for two compact metric spaces $X$ and $Y$.
By Lemma \ref{lem:precsim-cpt}, there are $1$-Lipschitz
surjective maps $\varphi : F \to Y$ and $\psi : G \to X$
such that $F$ is a closed subset of $X$
and that $G$ is a closed subset of $Y$.
Setting $X' := \varphi^{-1}(G)$, we see that
$f := \psi\circ\varphi : X' \to X$ is a $1$-Lipschitz surjective map.
A standard discussion (see \cite{BBI}*{Section 1.6})
proves that $X' = X$ and $f$ is an isometry.
Therefore, $F = X$, $G = Y$, and $\varphi$, $\psi$ are isometries.
This completes the proof.
\end{proof}

\begin{defn}[$\varepsilon$-Wider space]
Let $X$ and $Y$ be two extended metric spaces,
and let $\varepsilon \ge 0$.
We say that $Y$ is \emph{$\varepsilon$-wider than $X$}
if there exists a map $f : X \to Y$ such that
\begin{equation}
\label{eq:wider}
d_X(x,x') \le d_Y(f(x),f(x')) + \varepsilon
\end{equation}
for any $x,x' \in X$.
We call $f(X)$ the \emph{target} of $X$.
If $Y$ is $0$-wider than $X$, then we say that
$Y$ is \emph{wider than $X$}.
\end{defn}

If $X$ is $\varepsilon_1$-wider than $Y$
and if $Y$ is $\varepsilon_2$-wider than $Z$,
then $X$ is $(\varepsilon_1+\varepsilon_2)$-wider than $Z$.
If $\dGH(X,Y) < \varepsilon$ for two compact metric spaces
$X$ and $Y$, then
$X$ is $2\varepsilon$-wider than $Y$.

\begin{rem}
Consider the following conditions (i) and (ii).
\begin{enumerate}
\item[(i)] $X \precsim Y$.
\item[(ii)] For any $\varepsilon > 0$ and for any finite subspace
$X' \subset X$, the $Y$ is $\varepsilon$-wider than $X'$.
\end{enumerate}
In the case where $X$ is a (non-extended) metric space,
(i) and (ii) are equivalent to each other.
In the case where $X$ is extended,
we have (ii)$\Rightarrow$(i),
however (i)$\Rightarrow$(ii) does not necessarily hold
(because \eqref{eq:wider} does not hold if the left-hand side is infinity
and if the right-hand side is finite).
\end{rem}

Let $\cH$ denote the GH space, and
$\cH(N)$ the set of isometry classes of metric spaces consisting of
at most $N$ points.
We see that $\bigcup_{N=1}^\infty \cH(N)$ is a dense subset of $\cH$.
For a real number $D \ge 0$ and a metric space $X$,
we define
$d_{X\wedge D}(x,x') := d_X(x,x') \wedge D$ for $x,x' \in X$,
and $X \wedge D := (X,d_{X \wedge D})$.

\begin{lem} \label{lem:GH-precsim}
Let $X$ and $Y$ be compact metric spaces
such that $\dGH(X,Y) < \varepsilon$
for a real number $\varepsilon > 0$.
Then, for any compact metric space $Y'$ with $Y' \precsim Y$,
there exists a compact metric space $X'$ such that
\begin{enumerate}
\item $\diam X' \le \diam Y'$,
\item $X' \precsim X$,
\item $\dGH(X',Y') \le 3\varepsilon$,
\item if $Y' \in \cH(N)$, then  we have $X' \in \cH(N)$.
\end{enumerate}
\end{lem}

\begin{proof}
By the compactness of $Y'$, there is a finite $\varepsilon$-net $Y''$ of $Y'$,
i.e, the $\varepsilon$-neighborhood of $Y''$ covers $Y'$.
If $Y'$ itself is finite, then we assume $Y'' = Y'$.
The $X$ is $2\varepsilon$-wider than $Y''$ and we denote
the target of $Y''$ by $X''$.
Put $N := \# Y''$ and embed
$Y''$ into $(\R^N,\|\cdot\|_\infty)$
by the Kuratowski embedding (see \cite{Hn:G-emb}).
Then we have a map
$f : X'' \to Y'' \hookrightarrow (\R^N,\|\cdot\|_\infty)$
that is $1$-Lipschitz up to $2\varepsilon$.
We apply Lemma \ref{lem:f-tilde} to obtain
a $1$-Lipschitz map $\tilde{f} : X'' \to (\R^N,\|\cdot\|_\infty)$
such that $\|\tilde{f}(x)-f(x)\|_\infty \le 2\varepsilon$
for any $x,x' \in X''$.
It then holds that
$\tilde{f}(X'') \precsim X$ and
$\dH(\tilde{f}(X''),Y') \le 3\varepsilon$.
The space $X' := \tilde{f}(X'') \wedge \diam Y'$
is required one.
This completes the proof.
\end{proof}

\begin{defn}[Pyramid] \label{defn:pyramid}
A subset $\cP \subset \cH$ is called a \emph{pyramid}
if it satisfies the following (i)--(iii).
\begin{enumerate}
\item[(i)] If $X \precsim Y \in \cP$, then $X \in \cP$.
\item[(ii)] For any $\varepsilon > 0$ and $X, Y \in \cP$,
there exist $X', Y', Z \in \cP$ such that
$\dGH(X,X') \le \varepsilon$,
$\dGH(Y,Y') \le \varepsilon$,
and $X', Y' \precsim Z$.
\item[(iii)] $\cP$ is a nonempty GH closed subset of $\cH$.
\end{enumerate}
Denote by $\Pi$ the set of pyramids.
\end{defn}

Gromov's definition of pyramid for metric measure spaces
\cite{Gmv:green}
is similar to Definition \ref{defn:pyramid}.
However, the order relation $\precsim$
and Definition \ref{defn:pyramid}(ii) are different from Gromov's.
Also, (iii) is missing in Gromov's definition
(in fact, later added in our previous work \cites{Sy:mmg, Sy:mmlim})

For an extended metric space $X$, we define
\[
\cP_X := \{\,Y \in \cH \mid Y \precsim X\,\}.
\]
The following lemma is easy to prove.

\begin{lem} \label{lem:procsim-lim}
Let $X, Y, X_n, Y_n$, $n=1,2,\dots$, be compact metric spaces
such that $X_n$ and $Y_n$ GH converge to $X$ and $Y$
as $n\to\infty$, respectively.
If $X_n \precsim Y_n$ for any $n$, then $X \precsim Y$.
\end{lem}

We prove the following.

\begin{prop}
For any extended metric space $X$,
the set $\cP_X$ is a pyramid.
\end{prop}

\begin{proof}
Let $S$ be an extended metric space.

It is clear that $\cP_S$ satisfies Definition \ref{defn:pyramid}(i).

We prove (ii) for $\cP_S$.
Note that this is nontrivial if $S$ is noncompact or extended.
Take any $\varepsilon > 0$, $X, Y \in \cP_S$.
Put $D := \diam X \vee \diam Y$.
Let $X'$ (resp.~$Y'$) be a finite $\varepsilon$-net of $X$ (resp.~$Y$).
It follows from $X, Y \precsim S$
that $S$ is $\varepsilon_n$-wider than $X'$ and $Y'$
for a sequence $\varepsilon_n \to 0+$.
Let $A_n, B_n \subset S$ be their targets.
Setting $Z_n := (A_n \cup B_n) \wedge D$, we see
that $Z_n \precsim S$ and $Z_n$ is $\varepsilon_n$-wider than 
$X'$ and $Y'$.
Since the number of points of $Z_n$ is bounded above,
by taking a subsequence,
$Z_n$ GH converges to a finite metric space, say $Z$.
Since $X', Y' \precsim Z \precsim S$,
we obtain (ii).

(iii) is clear.

This completes the proof.
\end{proof}

The following proposition directly follows from the transitivity of $\precsim$.

\begin{prop} \label{prop:PXsim}
For two extended metric spaces $X$ and $Y$,
we have
\begin{enumerate}
\item $X \precsim Y$ if and only if $\cP_X \subset \cP_Y$;
\item $X \sim Y$ if and only if $\cP_X = \cP_Y$. 
\end{enumerate}
\end{prop}

\begin{lem} \label{lem:pyramid-lim}
If a sequence of pyramids converges in the weak Hausdorff sense,
then the limit is a pyramid.
The set of pyramids is sequentially compact
with respect to the weak Hausdorff convergence. 
\end{lem}

\begin{proof}
Let us prove the first part of the lemma.
Assume that a sequence of pyramids $\cP_n$, $n=1,2,\dots$,
converges to a subset $\cP \subset \cH$ in the weak Hausdorff sense.

Definition \ref{defn:pyramid}(i) for $\cP$ follows from
Lemma \ref{lem:GH-precsim}.

We prove (ii) for $\cP$.
Take any $X, Y \in \cP$ and $\varepsilon > 0$ and fix them.
Let $X_\varepsilon$ (resp.~$Y_\varepsilon$) be
a finite $\varepsilon$-net of $X$ (resp.~$Y$).
There are $X_n, Y_n \in \cP_n$, $n=1,2,\dots$, such that
$\dGH(X_n,X_\varepsilon),\dGH(Y_n,Y_\varepsilon) < \varepsilon_n$
for a sequence $\varepsilon_n \to 0+$.
Condition (ii) for $\cP_n$ implies that
there are $X_n', Y_n', Z_n \in \cP_n$ such that
$\dGH(X_n,X_n') < \varepsilon_n$,
$\dGH(Y_n,Y_n') < \varepsilon_n$,
and $X_n', Y_n' \precsim Z_n$.
There are nets $X_{n,\varepsilon}', Y_{n,\varepsilon}'$
of $X_n', Y_n'$ respectively
such that $\dGH(X_{n,\varepsilon}',X_\varepsilon),
\dGH(Y_{n,\varepsilon}',Y_\varepsilon) < 3\varepsilon_n$,
$\# X_{n,\varepsilon}' = \# X_\varepsilon$, and
$\# Y_{n,\varepsilon}' = \# Y_\varepsilon$
for every sufficiently large $n$,
where $\# A$ denotes the number of elements of a set $A$.
It follows from $X_n', Y_n' \precsim Z_n$ that
$Z_n$ is wider than $X_{n,\varepsilon}'$ and $Y_{n,\varepsilon}'$.
Let $A_n,B_n \subset Z_n$ be their targets.
The union $A_n \cup B_n$ is a subset of $Z_n$
consisting of at most $(\# X_\varepsilon + \# Y_\varepsilon)$ points.
Let $D := \diam X_\varepsilon \vee \diam Y_\varepsilon$.
By taking a subsequence, $(A_n \cup B_n) \wedge D$ \; GH converges to
a finite metric space, say $Z$.
Since $X_\varepsilon, Y_\varepsilon \precsim Z$
(see Lemma \ref{lem:procsim-lim}),
we obtain (ii) for $\cP$.

(iii) for $\cP$ is clear.

We have proved that $\cP$ is a pyramid.

The later part follows from the first and
Lemma \ref{lem:weakH-cpt}.
\end{proof}

For an integer $N \ge 1$ and a real number $D \ge 0$,
we define
\[
\cH(N,D) := \{\,X \in \cH(N) \mid \diam X \le D\,\}.
\]
It is well-known that $\cH(N,D)$ is GH compact.

\begin{lem} \label{lem:cPn-conv-1}
Let $\cP$ and $\cP_n$, $n=1,2,\dots$, be pyramids.
Then, the following {\rm(1)} and {\rm(2)} are equivalent to each other.
\begin{enumerate}
\item $\cP_n$ converges to $\cP$ in the weak Hausdorff sense
as $n\to\infty$.
\item For any $N \ge 1$ and $D \ge 0$,
the set $\cP_n \cap \cH(N,D)$ converges to $\cP \cap \cH(N,D)$
in the Hausdorff sense as $n\to\infty$.
\end{enumerate} 
\end{lem}

\begin{proof}
We prove (1)$\Rightarrow$(2).
Assume (1).
We take any $N \ge 1$, $D \ge 0$, and fix them.
Since $\cP_n \cap \cH(N,D)$ and $\cP \cap \cH(N,D)$ are compact,
it suffices to prove the weak Hausdorff convergence
of $\cP_n \cap \cH(N,D)$ to $\cP \cap \cH(N,D)$
(see Lemma \ref{lem:weak-H}).

By (1), for any $X \in \cP \cap \cH(N,D)$,
there is a sequence $X_n \in \cP_n$, $n=1,2,\dots$,
GH converging to $X$.
It follows from Lemma \ref{lem:GH-precsim}
that there is a sequence $Y_n \in \cH(N)$ with $Y_n \precsim X_n$
such that $Y_n$ GH converges to $X$ and $\diam Y_n \le D$.
Since $Y_n \in \cP_n \cap \cH(N,D)$,
we have $\dGH(X,\cP_n \cap \cH(N,D)) \to 0$ as $n\to\infty$.

Take any $X \in \cH \setminus (\cP \cap \cH(N,D))$.
If $X$ does not belong to $\cH(N,D)$, then
$\liminf_{n\to\infty} \dGH(X,\cP_n \cap \cH(N,D))
\ge \dGH(X,\cH(N,D)) > 0$.
If $X$ belongs to $\cH(N,D)$, then $X$ does not belong to $\cP$,
which together with (1) implies\\
$\liminf_{n\to\infty} \dGH(X,\cP_n \cap \cH(N,D))
\ge \liminf_{n\to\infty} \dGH(X,\cP_n) > 0$.

We have proved that $\cP_n \cap \cH(N,D)$
converges to $\cP \cap \cH(N,D)$ in the weak Hausdorff sense.
This together with Lemma \ref{lem:weak-H} implies (2).

We prove (2)$\Rightarrow$(1).

Take any $X \in \cP$.
There is a sequence $X_N \in \cP \cap \cH(N,N)$, $N=1,2,\dots$,
GH converging to $X$.
By (2), for each $N$ there is a sequence
$X_{N,n} \in \cP_n \cap \cH(N,N)$, $n=1,2,\dots$,
GH converging to $X_N$.
There is a sequence $N(n)$ such that
$X_{N(n),n}$ GH converges to $X$ as $n\to\infty$.
We therefore have $\dGH(X,\cP_n) \to 0$ as $n\to\infty$.

Suppose that there is $X \in \cH \setminus \cP$ such that
$\liminf_{n\to\infty} \dGH(X,\cP_n) = 0$.
Taking a subsequence, we assume
$\lim_{n\to\infty} \dGH(X,\cP_n) = 0$.
There is a sequence $X_n \in \cP_n$, $n=1,2,\dots$,
GH converging to $X$.
Setting $\delta := \dGH(X,\cP)$, we have $\delta > 0$
because of the closedness of $\cP$.
Put $D := \diam X$.
There are $N$ and $X' \in \cH(N,D)$ such that
$X' \precsim X$ and $\dGH(X,X') < \delta$.
By applying Lemma \ref{lem:GH-precsim},
there is a sequence $X_n' \in \cH(N,D)$, $n=1,2,\dots$,
such that $X_n' \precsim X_n$
and $X_n'$ GH converges to $X'$.
It follows from (2) that $X' \in \cP$,
which contradicts $\dGH(X,X') < \delta$.

We have proved (1).
This completes the proof of the lemma.
\end{proof}

\begin{lem} \label{lem:cPn-conv-2}
Let $0 \le D \le D'$, $N \le N'$, and $\cP,\cP_n \subset \cH$,
where $n=1,2,\dots$.
Then, if $\cP_n \cap \cH(N',D')$ Hausdorff converges to
$\cP \cap \cH(N',D')$ as $n\to\infty$,
then
$\cP_n \cap \cH(N,D)$ Hausdorff converges to $\cP \cap \cH(N,D)$
as $n\to\infty$.
\end{lem}

\begin{proof}
Assume that $\cP_n \cap \cH(N',D')$ Hausdorff converges to
$\cP \cap \cH(N',D')$.
By the compactness of $\cH(N,D)$,
replacing $\{\cP_n\}$ with a subsequence
we see that $\cP_n \cap \cH(N,D)$ Hausdorff converges
to a subset of $\cH(N,D)$, say $\cP_\infty$.
It suffices to prove that $\cP_\infty = \cP \cap \cH(N,D)$.
By the assumption, we have $\cP_\infty \subset \cP$.
For any $X \in \cP \cap \cH(N,D)$,
since $X \in \cP \cap \cH(N',D')$,
there is a sequence $X_n \in \cP_n \cap \cH(N',D')$, $n=1,2,\dots$,
GH converging to $X$.
Applying Lemma \ref{lem:GH-precsim} yields that
there is a sequence $Y_n \in \cH(N,D)$, $n=1,2,\dots$, with
$Y_n \precsim X_n$ such that $Y_n$ GH converges to $X$.
Since $Y_n \in \cP_n \cap \cH(N,D)$, we have $X \in \cP_\infty$.
We thus obtain $\cP_\infty = \cP \cap \cH(N,D)$.
This completes the proof.
\end{proof}

\begin{defn}[Pseudometric $\rho$] \label{defn:rho}
For two extended metric spaces $X$ and $Y$, we define
\begin{align*}
\rho_N(X,Y) &:= \dH(\cP_X \cap \cH(N,N),\cP_Y \cap \cH(N,N)),\\
\rho(X,Y) &:= \sum_{N=1}^\infty 2^{-N} \rho_N(X,Y),
\end{align*}
where $\dH$ is the Hausdorff distance with respect to the GH metric.
\end{defn}


Since $\rho_N \le N$ we see $\rho \le 2$.
It is obvious that $\rho_N$ and $\rho$ are pseudometrics.

\begin{thm} \label{thm:rho-conv}
Let $X_n$ and $Y$ be extended metric spaces, $n=1,2,\dots$.
Then the following {\rm(1)} and {\rm(2)} are equivalent to each other.
\begin{enumerate}
\item $\cP_{X_n}$ converges to $\cP_Y$ in the weak Hausdorff sense
as $n\to\infty$.
\item $X_n$ $\rho$-converges to $X$ as $n\to\infty$.
\end{enumerate}
\end{thm}

\begin{proof}
It follows from Lemmas \ref{lem:cPn-conv-1} and \ref{lem:cPn-conv-2}
that
(1) is equivalent to the $\rho_N$-convergence of $X_n$ to $Y$
for any $N$, which is equivalent to (2).
This completes the proof.
\end{proof}

\begin{cor} \label{cor:rho-sim}
For two extended metric spaces $X$ and $Y$,
we have $\rho(X,Y) = 0$ if and only if $X \sim Y$.
\end{cor}

\begin{proof}
The corollary follows from
Proposition \ref{prop:PXsim} and Theorem \ref{thm:rho-conv}.
\end{proof}

Corollary \ref{cor:rho-sim} implies that
$\rho$ induces a metric on $\cM/\!\sim$,
which we denote by the same notation $\rho$.

\begin{lem} \label{lem:dH-cP-dGH}
Let $X$ and $Y$ be two compact metric spaces.
Then we have
\begin{align}
\tag{1} \dH(\cP_X,\cP_Y) &\le 3\dGH(X,Y),\\
\tag{2} \rho(X,Y) &\le 3\dGH(X,Y),
\end{align}
where $\dH$ is the Hausdorff distance with respect to $\dGH$.
\end{lem}

\begin{proof}
(1) follows from Lemma \ref{lem:GH-precsim}.

It follows from Lemma \ref{lem:GH-precsim}
that $\rho_N(X,Y) \le 3\dGH(X,Y)$, which implies (2).
\end{proof}

\begin{lem} \label{lem:emb}
Let $X_n$ and $Y$ be compact metric spaces, $n=1,2,\dots$.
Then the following {\rm(1)} and {\rm(2)} are equivalent to each other.
\begin{enumerate}
\item $X_n$ GH converges to $Y$ as $n\to\infty$.
\item $\cP_{X_n}$ converges to $\cP_Y$ in the weak Hausdorff sense
as $n\to\infty$.
\end{enumerate}
\end{lem}

\begin{proof}
(1)$\Rightarrow$(2) follows from Lemma \ref{lem:dH-cP-dGH}(1).

To prove (2)$\Rightarrow$(1), we assume (2).

We first prove the GH-precompactness of $\{X_n\}$.
If not, then, by taking a subsequence of $\{X_n\}$,
there is a sequence of $\delta$-discrete nets
$X_n' \subset X_n$, $n=1,2,\dots$, for a real number $\delta > 0$
such that $\# X_n' \to \infty$ as $n\to\infty$.
Take any $N \ge 1$ and $D \ge \delta$ and fix them.
Let $n$ be so large that $\# X_n' \ge N$.
We take a subset $X_n'' \subset X_n'$ with $\# X_n'' = N$.
We have $X_n''\wedge D \in \cP_{X_n} \cap \cH(N,D)$,
which together with Lemma \ref{lem:cPn-conv-1} implies
that there is $Y_n \in \cP_Y$ such that
$\dGH(X_n''\wedge D,Y_n) \to 0$.
Thereby, $Y$ is $\varepsilon_n$-wider than $X_n''\wedge D$
for a sequence $\varepsilon_n \to 0$
and so the target is $(\delta/2)$-discrete net with cardinality $N$
for every sufficiently large $n$.
By the arbitrariness of $N$,
this is a contradiction to the compactness of $Y$.
Therefore, $\{X_n\}$ is GH-precompact.

Suppose that $X_n$ does not GH converge to $Y$.
By taking a subsequence, $X_n$ GH converges to
a compact metric space, say $Z$, that is not isometric to $Y$.
Applying (1)$\Rightarrow$(2) yields that
$\cP_{X_n}$ converges to both $\cP_Y$ and $\cP_Z$
in the weak Hausdorff sense.
By Theorem \ref{thm:rho-conv},
$X_n$ $\rho$-converges to both $Y$ and $Z$,
which implies $\rho(Y,Z) = 0$.
By Corollary \ref{cor:rho-sim} and Lemma \ref{lem:prec-isom},
they are isometric to each other.  This is a contradiction.
The proof of the lemma is completed.
\end{proof}

\begin{lem} \label{lem:dense}
The image of the embedding map
$\cH \ni X \mapsto \cP_X \in \Pi$ is dense in $\Pi$.
\end{lem}

\begin{proof}
Let $\cP$ be a pyramid.
By the separability of $\cH$,
there is a dense countable subset $\{Y_n\} \subset \cP$.
Take any $\varepsilon > 0$ and fix it.
Let us define a sequence $Z_n \in \cP$, $n=1,2,\dots$, inductively.
We define $Z_1 := Y_1$.
If $Z_n$ is defined for a number $n$,
then, by the definition of the pyramid,
there are $Z_n', Y_{n+1}', Z_{n+1} \in \cP$ such that
$\dGH(Z_n,Z_n'), \dGH(Y_{n+1},Y_{n+1}') < 3^{-1}2^{-n}\varepsilon$ and
$Z_n', Y_{n+1}' \precsim Z_{n+1}$.

Now, let $n$ be any integral number with $n \ge 2$.
Since $Y_n' \in \cP_{Z_n}$,
we see $\dGH(Y_n,\cP_{Z_n}) \le 3^{-1}2^{-n+1}\varepsilon$.
By Lemma \ref{lem:GH-precsim}, we have
$\cP_{Z_n} \subset B_{2^{-n}\varepsilon}(\cP_{Z_n'})
\subset B_{2^{-n}\varepsilon}(\cP_{Z_{n+1}})$,
so that
$Y_n \in B_\varepsilon(\cP_{Z_m})$ for any $m \ge n$.
By taking a subsequence of $\{\cP_{Z_m}\}$,
it converges to a pyramid, say $\cP_\varepsilon$,
in the weak Hausdorff sense.

We prove $Y_n \in B_\varepsilon(\cP_\varepsilon)$ in the following.
Let $\delta$ be any positive real number.
There is $W_m \in \cP_{Z_m}$ such that $\dGH(Y_n,W_m) < \varepsilon + \delta$.
We take a $\delta$-net $\mathcal{N}$ of $Y_n$.
There is a net $W_m' \subset W_m$ such that $\# W_m' \le \#\mathcal{N}$
and $\dGH(\mathcal{N},W_m') < \varepsilon + \delta$.
Taking a subsequence of $\{W_m'\}$, it GH converges to
a finite metric space, say $W$, which belongs to $\cP_\varepsilon$.
Since $\dGH(Y_n,W) \le \varepsilon + 2\delta$,
we have $Y_n \in B_{\varepsilon+2\delta}(\cP_\varepsilon)$.
By the arbitrariness of $\delta$, we obtain
$Y_n \in B_\varepsilon(\cP_\varepsilon)$.

Since $\{Y_n\}$ is dense in $\cP$, we have
$\cP \subset B_\varepsilon(\cP_\varepsilon)$.
By $Z_n \in \cP$
we also have $\cP_\varepsilon \subset \cP$.
As $\varepsilon \to 0$,
the $\cP_\varepsilon$ converges to $\cP$.
This completes the proof.
\end{proof}

We have the following theorem.
Refer to \cite{KL}*{Section 3} for ultralimit.

\begin{thm} \label{thm:pyramid-ultralimit}
Let $\{X_n\}$ be a sequence of extended metric spaces
such that $\cP_{X_n}$ converges to a pyramid $\cP$
in the weak Hausdorff sense.
Then, any ultralimit $X_\infty$ of $\{X_n\}$ satisfies
$\cP = \cP_{X_\infty}$.
In particular, all ultralimits of $\{X_n\}$ are equivalent to each other.
\end{thm}

\begin{proof}
We prove $\cP \subset \cP_{X_\infty}$.
Take any $Y \in \cP$.
Since $\cP_{X_n} \to \cP$, there is a sequence
$Y_n \in \cP_{X_n}$, $n=1,2,\dots$,
GH converging to $Y$.
For any finite subset $Y' \subset Y$, we see that
$Y_n$ is $\varepsilon_n$-wider than $Y'$
for a sequence $\varepsilon_n \to 0+$,
so that $X_n$ is $\varepsilon_n$-wider than $Y'$,
which implies $Y' \precsim X_\infty$ and so $Y' \in \cP_{X_\infty}$.
Since $\cP_{X_\infty}$ is closed and $Y'$ approximates $Y$,
we have $Y \in \cP_{X_\infty}$.

We prove $\cP_{X_\infty} \subset \cP$.
Take any $Y \in \cP_{X_\infty}$ and put $D := \diam Y$.
We also take any finite subset $Y' \subset Y$ and put $N := \# Y'$.
The ultralimit $X_\infty$ is wider than $Y'$ and
we denote the target of $Y'$ by $X_\infty'$.
There is a sequence $X_n' \subset X_n$, $n=1,2,\dots$,
with $\# X_n' = N$
that is $\varepsilon_n$-wider than $X_\infty'$
for a sequence $\varepsilon_n \to 0+$.
Let $X_n'' := X_n'\wedge D$.
Taking a subsequence, $X_n''$ GH converges to
a finite metric space, say $X'$.
Since $X_n'' \in \cP_{X_n}$ we see $X' \in \cP$.
Moreover we have $Y' \precsim X'$
and therefore $Y' \in \cP$.
Since $Y'$ approximates $Y$, we obtain $Y \in \cP$.
This completes the proof.
\end{proof}

\begin{cor} \label{cor:PX}
For any pyramid $\cP$ there exists an extended metric space $X$
such that $\cP = \cP_X$.
\end{cor}

\begin{proof}
Let $\cP$ be any pyramid.
By Lemma \ref{lem:dense},
there is a sequence $X_n \in \cH$ such that
$\cP_{X_n}$ converges to $\cP$ in the weak Hausdorff sense.
Theorem \ref{thm:pyramid-ultralimit} completes the proof.
\end{proof}

Let $[X] \in \cM/\!\sim$ denote the equivalence class
represented by $X \in \cM$.

\begin{proof}[Proof of Main Theorem \ref{mthm:main}]
It follows from Proposition \ref{prop:PXsim}
and Corollary \ref{cor:PX} that
the map
\begin{equation} \label{eq:bij}
\cM/\!\sim \; \ni [X] \longmapsto \cP_X \in \Pi
\end{equation}
is bijective.
Combining Lemma \ref{lem:pyramid-lim} and
Theorem \ref{thm:rho-conv} yields that
$(\cM/\!\sim,\rho)$ is a compact metric space,
i.e., we obtain (1) of the theorem.

(2) follows from Lemmas \ref{lem:emb} and \ref{lem:dH-cP-dGH}(2).

(3) follows from Lemma \ref{lem:dense}.

(4) follows from Theorem \ref{thm:pyramid-ultralimit}.

This completes the proof.
\end{proof}

\begin{rem} \label{rem:cM}
We see that $\cM$ contains $(\cM/\!\sim,\rho)$,
which violates the axiom of regularity.
This means that neither $\cM$ nor $\cM/\!\sim$ are not sets.
To avoid it, we define $\cM/\!\sim$ to be $\Pi$
by the bijection \eqref{eq:bij}.
Then Main Theorem \ref{mthm:main} remains true rigorously in ZFC.
\end{rem}

\section{Examples} \label{sec:examples}

For $0 < D \le +\infty$,
we denote by $\Sigma_\infty(D)$
a countably infinite (extended) metric space such that
any two distinct points have distance $D$,
and by $\Sigma_n(D)$ a subspace of $\Sigma_\infty(D)$
consisting of exactly $n$ points.
For $0 < R < +\infty$ and $1 \le n \le \infty$,
we connect each point in $\Sigma_n(2R)$
to a fixed point, say $o$, 
(which is taken outside $\Sigma_n(2R)$) by a line segment
of length $R$ to obtain a geodesic metric space.
Such the geodesic metric space is called the \emph{$(n,R)$-spider}
and denoted by $Sp_n(R)$.

Note that $\Sigma_\infty(D)$, $0 < D \le +\infty$, is the maximum
with respect to $\precsim$
in the family of (extended) metric spaces with diameter $\le D$.
For any $0 < R < +\infty$,
the spider $Sp_\infty(R)$ is equivalent to $\Sigma_\infty(2R)$.
Note that $Sp_n(R)$ is not equivalent to $\Sigma_n(2R)$
for $n < +\infty$.
The following proposition characterizes maximal extended metric spaces.

\begin{prop} \label{prop:S-infty-infty}
An extended metric space $X$ is equivalent to $\Sigma_\infty(\infty)$
if and only if either or both of the following {\rm(1)} and {\rm(2)}
is satisfied.
\begin{enumerate}
\item $X$ has a subspace isometric to $\Sigma_\infty(\infty)$.
\item $X$ has a subspace isometric to
a {\rm(}non-extended{\rm)} metric space with infinite diameter.
\end{enumerate}
\end{prop}

\begin{proof}
If either or both of (1) and (2) is satisfied,
then it is easy to see that $X$ is equivalent to $\Sigma_\infty(\infty)$.

Assume that an extended metric space $X$ is equivalent
to $\Sigma_\infty(\infty)$.
We consider the equivalence relation on $X$
defined by $d_X(x,x') < \infty$ for $x,x' \in X$.
If we have infinitely many equivalence classes in $X$,
then (1) holds.
If the equivalence classes of $X$ are only finitely many,
then there is a equivalence class $X'$ of $X$
such that $\Sigma_\infty(\infty) \precsim X'$.
The diameter of $X'$ is infinity.
This completes the proof.
\end{proof}

Let $\{X_n\}_{n=1}^\infty$ be a sequence of extended metric spaces
satisfying that $X_n \precsim X_{n+1}$ for any $n$.
Then, since $\cP_{X_n} \subset \cP_{X_{n+1}}$,
the pyramid $\cP_{X_n}$ converges to the GH closure of
$\bigcup_{n=1}^\infty \cP_{X_n}$ in the weak Hausdorff sense.

We consider the particular case where
each $X_n$ is a subspace of a fixed extended metric space $X$
and satisfies $X_n \subset X_{n+1}$ for every $n$
and where $\bigcup_{n=1}^\infty X_n$ is a dense subset of $X$.
In this case, $\cP_{X_n}$ converges to $\cP_X$
and so $X_n$ $\rho$-converges to $X$ as $n\to\infty$.
The following are such specific examples.

\begin{ex} \label{ex:rho-conv}
\begin{enumerate}
\item \label{it:spider} For any $0 < R < +\infty$,
the spider $Sp_n(R)$ $\rho$-converges to
$Sp_\infty(R) \sim \Sigma_\infty(2R)$
as $n\to\infty$.
\item Let $S^n(1)$ denote a unit sphere in the Euclidean space
$\R^{n+1}$ and $S^\infty(1)$ a unit sphere
in an infinite-dimensional separable Hilbert space.
We equip them to the geodesic distance functions.
Then, $S^n(1)$ $\rho$-converges to $S^\infty(1)$.
\item Identifying antipodal pairs of points in the unit spheres,
we obtain real projective spaces $\RP^n$, $1 \le n \le \infty$.
We see that $\RP^n$ $\rho$-converges to $\RP^\infty$,
while $\RP^\infty$ is equivalent to $\Sigma_\infty(\pi/2)$.
\item The $n$-dimensional simplex $\rho$-converges
to the infinite-dimensional simplex as $n\to\infty$,
where we define metrics on simplexes in such a way that
the distance between any two distinct vertices is equal to one
and that each face has flat metric.
The infinite-dimensional simplex is equivalent to $\Sigma_\infty(1)$.
\item Let $X$ be a compact metric space consisting of
at least two distinct points, and let $1 \le p \le +\infty$.
Denote by $X_p^n$ the $n$-th product space $X^n$
equipped with the $l^p$-metric.
Then, $X_p^n$ $\rho$-converges to $X_p^\infty$ as $n\to\infty$.
If $p < +\infty$, then Proposition \ref{prop:S-infty-infty} implies
that $X_p^\infty \sim \Sigma_\infty(\infty)$.
If $p = +\infty$, then $X_p^\infty \sim \Sigma_\infty(D)$, $D := \diam X$,
because $X_p^\infty$ contains a subspace isometric to $\Sigma_\infty(D)$.
\end{enumerate}
\end{ex}

\begin{rem} \label{rem:dGH-rho-1}
We see that $\dGH(Sp_m(R),Sp_n(R)) = R$ for $m \neq n$.
This gives an example of a pair of $X$ and $Y$
such that $\dGH(X,Y)$ is not small,
but $\rho(X,Y)$ is arbitrarily small.
In particular, we do not have $\dGH(X,Y) \le C\rho(X,Y)$
in general, where $C$ is a constant.
\end{rem}

\begin{rem} \label{rem:dGH-rho-2}
More generally, if you have a compactification $\mathcal{C}$
of a noncompact complete metric space $(M,d_M)$,
then any metric $d_{\mathcal{C}}$ on $\mathcal{C}$
does not satisfy $d_M \le f \circ d_{\mathcal{C}}$,
where $f : [\,0,+\infty\,) \to \R$ is a function with $f(0+) = 0$.
In fact, there is a sequence of points in $M$
converging to a point in $\mathcal{C} \setminus M$.
Such a sequence is $d_{\mathcal{C}}$-Cauchy,
but not $d_M$-Cauchy,
so that $d_M \le f \circ d_{\mathcal{C}}$ does not hold.
\end{rem}

\section{Construction of metric $\rho_0$}
\label{sec:sigma}

In this section, we define a pointed version of the metric $\rho$,
which is applied to define the metric $\rho_0$ on
the family of pointed geodesic metric spaces.
We also obtain some quantitative estimate on
the $\rho$ and $\rho_0$ distances.

Denote by $\cM_0$ the family of isometry classes of
pointed extended metric spaces,
and by $\cH_0$ the family of isometry classes of
pointed compact metric spaces.

\begin{defn}[Pointed Gromov-Hausdorff distance \cite{Gmv:poly}]
The \emph{pointed GH distance} $\dGH((X,x_0),(Y,y_0))$
between two pointed compact metric spaces $(X,x_0)$ and $(Y,y_0)$
is defined to be the infimum of
\[
\dH(\varphi(X),\psi(Y)) \vee d_Z(\varphi(x_0),\psi(y_0)),
\]
where $\varphi : X \to Z$ and $\psi : Y \to Z$ both run
over all isometric embeddings into compact metric spaces $Z$.
\end{defn}

The definition implies
\begin{equation} \label{eq:dGH-dpGH}
\dGH(X,Y) \le \dGH((X,x_0),(Y,y_0)).
\end{equation}

\begin{defn}[Pointed version of $\precsim$ and $\sim$] \label{defn:pointed-sim}
For two pointed extended metric spaces $(X,x_0)$ and $(Y,y_0)$,
we define $(X,x_0) \precsim (Y,y_0)$ by the condition that
for any $x_1,\dots,x_N \in X$
there are $y_{i,k} \in Y$ with $y_{0,k} = y_0$ such that
\[
d_X(x_i,x_j) \le \liminf_{k\to\infty} d_Y(y_{i,k},y_{j,k})
\]
for any distinct $i,j = 0,1,\dots,N$.
We say $(X,x_0) \sim (Y,y_0)$ if
$(X,x_0) \precsim (Y,y_0)$ and $(Y,y_0) \precsim (X,x_0)$ both hold.
\end{defn}

\begin{defn}[Pointed pyramid]
A subset $\cP \subset \cH_0$ is called a \emph{pointed pyramid}
if the following {\rm(i)--(iii)} are satisfied.
\begin{enumerate}
\item[(i)] If $(X,x_0) \precsim (Y,y_0) \in \cP$, then $(X,x_0) \in \cP$.
\item[(ii)] For any $\varepsilon > 0$, $(X,x_0), (Y,y_0) \in \cP$,
there exist $(X',x_0'), (Y',y_0'), (Z,z_0) \in \cP$ such that
$\dGH((X,x_0),(X',x_0')) \le \varepsilon$,
$\dGH((Y,y_0),(Y',y_0')) \le \varepsilon$,
and $(X',x_0'), (Y',y_0') \precsim (Z,z_0)$.
\item[(iii)] $\cP$ is a nonempty GH closed subset of $\cH_0$.
\end{enumerate}
\end{defn}

For an extended pointed metric space $(X,x_0)$,
\[
\cP_{(X,x_0)} := \{\,(Y,y_0) \in \cH_0 \mid (Y,y_0) \precsim (X,x_0)\,\}
\]
is a pyramid.
Let $\cH_0(N,D) := \{\,(X,x_0) \in \cH_0 \mid X \in \cH(N,D)\,\}$
for $N \ge 1$, $D \ge 0$.

\begin{defn}[Metric $\rho$]
For a positive integer $N$
and for two pointed extended metric spaces $(X,x_0)$ and $(Y,y_0)$,
we define
\begin{align*}
\rho_N((X,x_0),(Y,y_0)) &:= \dH(\cP_{(X,x_0)} \cap \cH_0(N,N),\cP_{(Y,y_0)} \cap \cH_0(N,N)),\\
\rho((X,x_0),(Y,y_0))
&:= \sum_{N=1}^\infty 2^{-N} \rho_N((X,x_0),(Y,y_0)),
\end{align*}
where $\dH$ is the Hausdorff distance with respect to
the pointed GH metric.
\end{defn}

In the same way as in Section \ref{sec:rho}, we observe that
$\rho$ is a pseudometric on $\cM_0$
with $\rho((X,x_0),(Y,y_0)) \le 2$.
We also see that
$(X,x_0) \sim (Y,y_0)$ $\Leftrightarrow$ $\rho((X,x_0),(Y,y_0)) = 0$,
so that $\rho$ induces a metric on $\cM_0/\!\sim$.
It follows from \eqref{eq:dGH-dpGH} that
\begin{equation}
\rho(X,Y) \le \rho((X,x_0),(Y,y_0)).
\end{equation}
We have the following theorem exactly in the same way
as in the proof of Main Theorem \ref{mthm:main}.

\begin{thm} \label{thm:pointed}
\begin{enumerate}
\item $(\cM_0/\!\sim,\rho)$ is a compact metric space.
\item The natural projection from the pointed GH space $\cH_0$
to $\cM_0/\!\sim$ is $3$-Lipschitz topological embedding map.
\item The image of the projection in {\rm(2)} is dense
in $(\cM_0/\!\sim,\rho)$.
\item[(4)] For any $\rho$-convergent sequence in $\cM_0$,
any pointed ultralimit of it is equivalent to the $\rho$-limit.
\end{enumerate}
\end{thm}

The following lemmas are needed for the definition
of $\rho_0$.

\begin{lem} \label{lem:conti}
Let $(X,x_0)$ be a pointed geodesic metric space.
Then, the pointed metric ball $(B(x_0,r),x_0)$
is $3$-Lipschitz continuous in $r > 0$ with respect to $\rho$.
\end{lem}

\begin{proof}
We estimate $\rho((B(x_0,r),x_0),(B(x_0,r+\varepsilon),x_0))$
for $r,\varepsilon > 0$.\\
Since $\cP_{(B(x_0,r),x_0)} \subset  \cP_{(B(x_0,r+\varepsilon),x_0)}$,
it suffices to prove that
\begin{equation} \label{eq:conti}
\cP_{(B(x_0,r+\varepsilon),x_0)} \cap \cH_0(N,N)
\subset B(\cP_{(B(x_0,r),x_0)} \cap \cH_0(N,N),3\varepsilon).
\end{equation}
Take any $(X',x_0') \in \cP_{(B(x_0,r+\varepsilon),x_0)} \cap \cH_0(N,N)$.
Then, $(B(x_0,r+\varepsilon),x_0)$ is wider than $(X',x_0')$.
Let $(\hat{X},x_0) \subset (B(x_0,r+\varepsilon),x_0)$ be
the target of $(X',x_0')$.
For each $x \in \hat{X}$ there is a point $y_x \in B(x_0,r)$
in a minimal geodesic joining $x_0$ and $x$
such that $d_X(x,y_x) \le \varepsilon$.
We put $\hat{Y} := \{\,y_x \mid x \in \hat{X}\,\}$.
It holds that $\hat{Y} \subset B(x_0,r)$ and
$\dGH((\hat{X},x_0),(\hat{Y},x_0)) \le \varepsilon$.
Applying Lemma \ref{lem:GH-precsim} yields that
there is a pointed finite metric space $(Y',y_0')$ satisfying
$(Y',y_0') \precsim (\hat{Y},x_0)$ and
$\dGH((X',x_0'),(Y',y_0')) \le 3\varepsilon$.
Since $(Y',y_0') \in \cP_{(B(x_0,r),x_0)} \cap \cH_0(N,N)$,
we obtain \eqref{eq:conti}.
This completes the proof.
\end{proof}

For $t > 0$ and a metric space $X$ with metric $d_X$,
we define $tX$ to be $X$ equipped with the metric $d_{tX} := t \, d_X$.

\begin{lem} \label{lem:scale}
Let $(X,x_0)$ be a pointed metric space with
$\diam X \le D$ for a constant $D$.
Then we have
\[
\rho((sX,x_0),(tX,x_0)) \le \frac12 D |s-t|
\]
for any $s,t > 0$.
\end{lem}

\begin{proof}
Take any $(X',x_0') \in \cP_{(X,x_0)}$.
Since $|s\,d_{X'}(x_1,x_2) - t\,d_{X'}(x_1,x_2)|
\le D|s-t|$ for any $x_1,x_2 \in X'$,
we have $\dGH((sX',x_0'),(tX',x_0')) \le D|s-t|/2$
(see Theorem \cite{BBI}*{7.3.25}).
This proves the lemma.
\end{proof}

\begin{lem} \label{lem:loc-Lip}
Let $(X,x_0)$ and $(Y,y_0)$ be two pointed geodesic metric spaces
and let $X(r) := (r^{-1}B(x_0,r),x_0)$ and $Y(r) := (r^{-1}B(y_0,r),y_0)$
for $r > 0$.
Then we have
\[
|\, \rho(X(r_1),Y(r_1)) - \rho((X(r_2),Y(r_2)) \,|
\le 8 \left| \, 1 - \frac{r_2}{r_1} \, \right|
\]
for any $r_1, r_2 > 0$.
\end{lem}

\begin{proof}
By Lemmas \ref{lem:conti} and \ref{lem:scale} together with
a triangle inequality,
\begin{align*}
\rho(X(r_1),X(r_2))
&\le \rho((r_1^{-1}B(x_0,r_1),x_0),(r_1^{-1}B(x_0,r_2),x_0))\\
&\ + \rho((r_1^{-1}B(x_0,r_2),x_0),(r_2^{-1}B(x_0,r_2),x_0))\\
&\le \frac{4}{r_1}|r_1 - r_2|,
\end{align*}
which implies
\begin{align*}
&|\, \rho(X(r_1),Y(r_1)) - \rho((X(r_2),Y(r_2)) \,|
\le \rho(X(r_1),X(r_2)) + \rho(Y(r_1),Y(r_2)) \\
&\le \frac{8}{r_1}|r_1 - r_2|.
\end{align*}
This completes the proof.
\end{proof}

\begin{defn}[Metric $\rho_0$]
For two pointed geodesic metric spaces $(X,x_0)$ and $(Y,y_0)$,
we define
\[
\rho_0((X,x_0),(Y,y_0))
:= \int_{(\,0,+\infty\,)} \rho((r^{-1}B(x_0,r),x_0),(r^{-1}B(y_0,r),y_0))
\, r e^{-r^2} \,dr.
\]
\end{defn}

Note that the integrability of the above
follows from Lemma \ref{lem:loc-Lip} and $\rho \le 2$.

\begin{defn}[Strongly equivalent]
For two pointed geodesic metric spaces $(X,x_0)$ and $(Y,y_0)$,
we say that $(X,x_0)$ is \emph{strongly equivalent} to $(Y,y_0)$
and write $(X,x_0) \simeq (Y,y_0)$ if
$(B(x_0,r),x_0) \sim (B(y_0,r),y_0)$ for any $r > 0$.
\end{defn}

Note that, in the case where $(X,x_0)$ and $(Y,y_0)$
are proper geodesic metric spaces,
we have $(X,x_0) \simeq (Y,y_0)$ if and only if they are
isometric to each other.

$\cG_0$ denotes the family of isometry classes
of pointed geodesic metric spaces.
Theorem \ref{thm:pointed} implies the following.

\begin{thm} \label{thm:sigma}
$\rho_0$ is a pseudometric on $\cG_0$.
We have $\rho_0((X,x_0),(Y,y_0)) = 0$ if and only if
$(X,x_0) \simeq (Y,y_0)$.
\end{thm}

We have the following.

\begin{thm} \label{thm:sigma-cpt}
$(\cG_0/\!\simeq,\rho_0)$ is a compact metric space.
For any $\rho_0$-convergent sequence of pointed geodesic
metric spaces, any pointed ultralimit of it is strongly equivalent to
the $\rho_0$-limit.
\end{thm}

\begin{proof}
Let $\{(X_n,x_{0,n})\}_{n=1}^\infty$ be a sequence
of pointed geodesic metric spaces.
Since $\{r^{-1}B(x_{0,n},r)\}_{n=1}^\infty$ has a $\rho$-convergent
subsequence for any $r > 0$, a diagonal argument yields
the existence of a subsequence $\{(X_{n(i)},x_{0,n(i)})\}_{i=1}^\infty$
such that $\{(r^{-1}B(x_{0,n(i)},r),x_{0,n(i)})\}_{i=1}^\infty$
$\rho$-converges for any rational number $r > 0$.
Let $(X_\infty,x_{0,\infty})$ be the pointed ultralimit
of $\{(X_{n(i)},x_{0,n(i)})\}_{i=1}^\infty$ with respect to
a non-principal ultrafilter $\omega$ on the set of positive integers.
Since $X_{n(i)}$ is a geodesic metric space, so is $X_\infty$
(see \cite{KL}*{Proposition 3.4}).
Therefore, the $\omega$-ultralimit of $(r^{-1}B(x_{0,n(i)},r),x_{0,n(i)})$
is isometric to $(r^{-1}B(x_{0,\infty},r),x_{0,\infty})$.
By Lemma \ref{lem:loc-Lip}, for any real number $r > 0$,
$(r^{-1}B(x_{0,n(i)},r),x_{0,n(i)})$ $\rho$-converges to
$(r^{-1}B(x_{0,\infty},r),x_{0,\infty})$ as $i\to\infty$.
Thus,  $(X_n,x_{0,n(i)})$ $\rho_0$-converges to
$(X_\infty,x_{0,\infty})$.
We obtain the compactness of $(\cG_0/\!\simeq,\rho_0)$.

The compatibility between the $\rho_0$-limit and the ultralimit
also follows from the above discussion.
This completes the proof.
\end{proof}

Denote by $\cP\cG_0$ the set of isometry classes
of pointed proper geodesic metric spaces.

\begin{thm} \label{thm:sigma-proper}
Let $(X_n,x_{0,n})$ and $(Y,y_0)$ be pointed proper geodesic
metric spaces, $n=1,2,\dots$.
Then the following {\rm(1)} and {\rm(2)} are equivalent to each other.
\begin{enumerate}
\item $(X_n,x_{0,n})$ pointed GH converges to $(Y,y_0)$
as $n\to\infty$.
\item $(X_n,x_{0,n})$ $\rho_0$-converges to $(Y,y_0)$
as $n\to\infty$.
\end{enumerate}
In other words, the quotient map $\iota : \cP\cG_0 \to \cG_0/\!\simeq$
is a topological embedding map.
\end{thm}

\begin{proof}
We prove (1)$\Rightarrow$(2).
Suppose that $(X_n,x_{0,n})$ pointed GH converges to
$(Y,y_0)$ as $n\to\infty$, but it does not $\rho_0$-converge
to $(Y,y_0)$.
By taking a subsequence, $(X_n,x_{0,n})$ $\rho_0$-converges to
a pointed metric space, say $(Z,z_0)$,
that is not strongly equivalent to $(Y,y_0)$.
By the assumption,
$(Y,y_0)$ is isometric to a ultralimit of $(X_n,x_{0,n})$,
which together with Theorem \ref{thm:sigma-cpt}
implies $(Y,y_0) \simeq (Z,z_0)$.
This is a contradiction.

We prove (2)$\Rightarrow$(1).
Assume that $(X_n,x_{0,n})$ $\rho_0$-converges to $(Y,y_0)$.
Then, $(r^{-1}B(x_{0,n},r),x_{0,n})$ $\rho$-converges to
$(r^{-1}B(y_0,r),y_0)$ for any $r > 0$.
Theorem \ref{thm:pointed}(2) proves that
$(r^{-1}B(x_{0,n},r),x_{0,n})$ pointed GH converges to
$(r^{-1}B(y_0,r),y_0)$ for any $r > 0$.
Therefore, $(X_n,x_{0,n})$ pointed GH converges to $(Y,y_0)$.
This completes the proof.
\end{proof}

\begin{proof}[Proof of Main Theorem \ref{mthm:main-pointed}]
(1) and (3) follow from Theorem \ref{thm:sigma-cpt}.
(2) follows from Theorem \ref{thm:sigma-proper}.
This completes the proof.
\end{proof}

\begin{ques}
Does the image of the embedding map
in Main Theorem \ref{mthm:main-pointed}
dense in $\cG_0/\!\simeq$\;?
\end{ques}

In the following,
we obtain some quantitative estimate of the $\rho$ and $\rho_0$
distances.

We call $\rad(X,x_0) := \sup_{x \in X} d_X(x_0,x)$
the \emph{radius} of a pointed metric space $(X,x_0)$.
It holds that $\diam X/2 \le \rad(X,x_0) \le \diam X$.

\begin{prop} \label{prop:rad}
Let $(X,x_0)$ be a pointed geodesic metric space
with $\rad(X,x_0) \ge R$ for a constant $R > 0$.
Then we have
\[
\rho(X,\Sigma_\infty(\infty))
\le \rho((X,x_0),(\Sigma_\infty(\infty),s_0)) \le (2+\sqrt{R})\,2^{1-\sqrt{R}},
\]
where $s_0 \in \Sigma_\infty(\infty)$.
\end{prop}

\begin{proof}
For any $\varepsilon$ with $0 < \varepsilon < R$,
there is a minimal geodesic $\gamma : [\,0,R-\varepsilon\,] \to X$
emanating from $x_0$ with length $R-\varepsilon$.
Put $N_0 := [(R-\varepsilon)^{1/2}]$.
For any $N$ with $1 \le N \le N_0$,
setting $x_i := \gamma(i(R-\varepsilon)/N)$, $i=1,\dots,N-1$,
we have $d_X(x_i,x_j) \ge (R-\varepsilon)/N \ge N$
for any distinct $i,j = 0,1,\dots,N-1$.
This proves that
$\cP_{(X,x_0)} \cap \cH_0(N,N) = \cP_{(\Sigma_N(N),s_N)} =
\cP_{(\Sigma_\infty(\infty),s_0)} \cap \cH_0(N,N)$,
where $s_N \in \Sigma_N(N)$ is a fixed point,
which implies $\rho_N((X,x_0),(\Sigma_\infty(\infty),s_0)) = 0$
for any $N$ with $1 \le N \le N_0$.
We also have $\rho_N \le N$ for all $N$.  Thus,
\[
\rho((X,x_0),(\Sigma_\infty(\infty),s_0))
\le \sum_{N=N_0+1}^\infty N 2^{-N} = (N_0+2)2^{-N_0}.
\]
By the arbitrariness of $\varepsilon$,
we obtain the proposition.
\end{proof}

\begin{lem} \label{lem:scale2}
Let $(X,x_0)$ and $(Y,y_0)$ be pointed metric spaces.
Then, for any $0 < r \le 1$ we have
\[
\rho((X,x_0),(Y,y_0)) \le \frac{1}{r}\, \rho((rX,x_0),(rY,y_0)).
\]
\end{lem}

\begin{proof}
Take any $\varepsilon$ with $\rho_N((rX,x_0),(rY,y_0)) < \varepsilon$.
Let $(X',x_0') \in \cP_{(X,x_0)} \cap \cH_0(N,N)$ be any pointed space.
Since $(rX',x_0') \in \cP_{(rX,x_0)} \cap \cH_0(N,N)$,
there is $(Y'',y_0'') \in \cP_{(rY,y_0)} \cap \cH_0(N,N)$ such that
$\dGH((rX',x_0'),(Y'',y_0'')) < \varepsilon$.
Setting $(Y',y_0') := ((r^{-1}Y'')\wedge N,y_0'')$,
we see that $(Y',y_0') \in \cP_{(Y,y_0)} \cap \cH_0(N,N)$ and
$\dGH((X',x_0'),(Y',y_0')) < \varepsilon/r$.
Since the same is true if we exchange $X$ and $Y$,
we have $\rho_N((X,x_0),(Y,y_0)) < \varepsilon/r$.
By the arbitrariness of $\varepsilon$,
we have
$\rho_N((X,x_0),(Y,y_0)) \le \rho_N((rX,x_0),(rY,y_0))/r$,
which implies the lemma.
\end{proof}

\begin{prop} \label{prop:rho-sigma}
Let $(X,x_0)$ and $(Y,y_0)$ be two pointed geodesic metric spaces
with $\rad(X,x_0), \rad(Y,y_0) \le R$ for a constant $R \ge 1/2$.
Then we have
\[
\rho(X,Y) \le \rho((X,x_0),(Y,y_0))
\le \frac{4R}{e^{-R^2}-e^{-4R^2}} \,\rho_0((X,x_0),(Y,y_0)).
\]
\end{prop}

\begin{proof}
For any $r$ with $r \ge R$ we have $B(x_0,r) = X$ and $B(y_0,r) = Y$,
which together with Lemma \ref{lem:scale2} implies
\begin{align*}
\rho_0((X,x_0),(Y,y_0))
&\ge \int_R^{2R} \rho((r^{-1}X,x_0),(r^{-1}Y,y_0)) \,r e^{-r^2} \, dr\\
&\ge \rho(((2R)^{-1}X,x_0),((2R)^{-1}Y,y_0))
\int_R^{2R} r e^{-r^2}\,dr\\
&\ge \frac{e^{-R^2}-e^{-4R^2}}{4R} \rho((X,x_0),(Y,y_0)).
\end{align*}
This completes the proof.
\end{proof}

\begin{rem}
\begin{enumerate}
\item The converse estimate of Proposition \ref{prop:rho-sigma}
does not hold.
In fact, if we take $s_0 \in Sp_\infty(\pi/2)$ as one of terminal points,
then we see $\rho((Sp_\infty(\pi/2),s_0),(\RP^\infty,x_0)) = 0$
for $x_0 \in \RP^\infty$.
However, we have
$\rho_0((Sp_\infty(\pi/2),s_0),(\RP^\infty,x_0)) > 0$.
\item In Proposition \ref{prop:rho-sigma}, the bound of the radius
is necessary.
In fact, denoting by $B^n(r)$ the $n$-dimensional Euclidean $r$-ball
centered at the origin $0$, we see that as $r\to +\infty$
the pointed ball $(B^n(r),0)$ converges
to $(\R^n,0)$ with respect to both $\rho$ and $\rho_0$.
On one hand, we have, as $r \to +\infty$,
\[
\rho((B^m(r),0),(B^n(r),0)) \to \rho((\R^m,0),(\R^n,0)) = 0
\]
because of $(\R^m,0) \sim (\R^n,0) \sim (\Sigma_\infty(\infty),x_0)$.
On the other hand, as $r \to +\infty$,
\[
\rho_0((B^m(r),0),(B^n(r),0)) \to \rho_0((\R^m,0),(\R^n,0)) > 0,
\]
because of $(\R^m,0) \not\simeq (\R^n,0)$ for $m\neq n$.
\end{enumerate}
\end{rem}

We next state a result on the stability of some properties
under convergence.
We say that $[X] \in \cM/\!\sim$ is \emph{geodesic}
if it contains a geodesic metric space.
We also say that $[X] \in \cM/\!\sim$ or $[(X,x_0)] \in \cG_0/\!\simeq$
is of \emph{Alexandrov curvature $\ge \kappa$}
(resp.~$\CAT(\kappa)$) for a real number $\kappa$
if it contains a space of Alexandrov curvature $\ge \kappa$
(resp.~$\CAT(\kappa)$ space).
We have the following.

\begin{cor}
Let $\kappa$ be a real number.
The properties of geodesic, of Alexandrov curvature $\ge \kappa$,
of $\CAT(\kappa)$ are respectively kept under
the $\rho$-convergence in $\cM$
and also under the $\rho_0$-convergence in $\cG_0$.
\end{cor}

\begin{proof}
It is well-known (and easy to prove) that the statement of the corollary
holds for ultralimits, which together with
Theorems \ref{mthm:main}(4) and \ref{mthm:main-pointed}(3)
implies the corollary.
\end{proof}

\end{document}